# PHYSICAL CONSTRAINTS OF NUMBERS

W. Mueckenheim

University of Applied Sciences, Augsburg

*All sciences need and many arts apply mathematics whereas mathematics seems to be independent of all of them, but only based upon logic. This conservative concept, however, needs to be revised because, contrary to Platonic idealism (frequently called "realism" by mathematicians), mathematical ideas, notions, and, in particular, numbers are not at all independent of physical laws and prerequisites.*

## ARE THERE ALL THE NATURAL NUMBERS?

Though this question seems to be too absurd to dwell upon, we will attempt it. The answer to the question obviously depends on the meaning of the notions "natural number" and "there is".

The *natural numbers* are defined by a set of axioms as given by G. Peano or by E. Schmidt, the most important ingredient of which is the principle of induction: If *n* is a natural number, then *n* + 1 is a natural number too. The smallest sequence

1, 2, 3, ...

which contains the unit 1 and obeys this principle is called "the sequence of natural numbers".

Set theory gives a meaning to "there is" by the axiom of infinity. *There is* a set *M* that contains the empty set { }, and that is such that if *m* belongs to *M*, then the union of *m* and {*m*} is also in *M*. Including one of the least natural of numbers, zero, the construction proposed by J. v. Neumann yields a definition of the set of all natural numbers

0 = { }, 1 = {{ }}, 2 = {{ },{{ }}}, ... .

The symbol of continuation appears here as above with different meaning though. Not the *potential* infinity of a *sequence* but the *actual* infinity of a *set* with cardinal number $\aleph_0$ is expressed by "...".

And God said: "Let there be light." And there was light. In this manner the axioms seem to create an infinitude of numbers. But there are set limits to profane creations. The *existence* of all natural numbers is not guaranteed by their *definition* alone. But in order to discuss this topic in detail we need some formalized concept of existence.

In its full generality the meaning of existence is a difficult philosophical question. The responses to it span the wide range between materialism and solipsism. But we need not consider this problem in great depth. What are the essentials of an idea? The answer to this question is comparatively easy. An idea exists on its own, if it is uniquely identifiable, i.e., if it can be distinguished from any other idea. A poem not yet written but existing in its author's mind is an idea as well as a mental image or a mathematical problem. The result of the latter is often a number. While numbers belong to the set of ideas, they have to satisfy even a stronger criterion in order to be considered as really existing. Like an ordinary idea a number *n* can be individualised by a mere name. But its existence is certainly not yet established by labels like "2" or "π". "Number" is a patent of nobility, not issued unless its *value*, i.e., its ratio with respect to the unit, *n*/1, *can be fixed precisely or at least to any desired precision.*

The basic and most secure method to establish the reality of a number is to form its fundamental set. The roman numerals are reminiscent of this method. While "2" is not a number but only a name, "II" is both, a name and a part of that number, namely of the fundamental set

2 = {all pairs like: II, you & me, mum & dad, sun & moon}.

2 is that property which all its subsets have in common. Of course this realization of 2 presupposes some *a priori* knowledge about 2. But here we are concerned with the mere realization. The same method can be applied in case of 3

3 = {all triples like: III, sun & moon & earth, father & son & holy spirit}.

Although IIII is easily intelligible, and pigeons are able to distinguish up to IIIIII slots at first glance, it is impractical to realize larger numbers in this way. The Romans ceased at IV already. And we would have great difficulties to identify numbers beyond 10 in v. Neumann's rather complex construction.

Position systems, decimal or binary or other *n*-adic systems, have the advantage to accomplish both identify a number and put it in order with other numbers by economical consumption of symbols.

Does the fundamental set of 4711 exist? We don't know. It did at least in Cologne at the beginning of the 19$^{th}$ century. A set of $10^{1000}$ elements does definitely not exist as will be shown below. Nevertheless 4711 and $10^{1000}$ are natural numbers. Their values, i.e., their ratios with respect to the unit are exactly determined.

It is impossible, however, to satisfy this condition for *all* natural numbers. It would require an unlimited amount of resources. But the universe is finite - at least that part available to us. Here is a simple means to realise the implications: First find out how many different natural numbers can be stored on a 10 GB hard disk. Then, step by step, expand the horizon to the $10^{11}$ neurons of your brain, to the $10^{28}$ atoms of your body, to the $10^{50}$ atoms of our earth, to the $10^{68}$ atoms of our galaxy, and finally to the $10^{78}$ protons within the universe. In principle, the whole universe could be turned into a big computer, but with far fewer resources than is usually expected without a thought be given to it.

Does a number exist in spite of the fact that we do not know and *cannot* know much more about it but that it should be a natural number? Does a poem exist, if nobody knows anything about it, except that it consists of 80 characters? Here is a poem which does exist:

>Flower, sweet primrose, show us your face;
>Now willow, t'is time to bud, make haste!

H. C. Andersen wrote it in his fairy tale "The Snowman". We can read, enjoy, love, learn, forget, discuss or criticize it. But that is not possible for the large majority of 80-letters strings, namely those which have not been written as yet.

The attempt to label each 80-digits number like the following

12345678901234567890123456789012345678901234567890123456789012345678901234567890

by a single proton only, would already consume more protons than the universe can supply. What you see is a number with no doubt. *Each* of those 80-digits strings can be noted on a small piece of paper. But *all* of them cannot even exist as *individual ideas* simultaneously, let alone as numbers with definite values. Given that photons and leptons can be recruited to store bits and given that the mysterious dark matter consists of particles which can be used for that purpose too, it is nevertheless quite impossible to encode the values of $10^{100}$ natural numbers in order to have them simultaneously available. (An advanced estimation, based on the Planck-length, leads to an upper limit of $10^{205}$ but the concrete number is quite irrelevant. So let us stay with $10^{100}$, the Googol.)

Even some single numbers smaller than $2^{10^{100}}$ can never be stored, known, or thought of. In short they do not exist. (It must not be forgotten: Also one's head, brain, mind, and all thoughts belong to the interior of the universe.) To avoid misunderstanding: Of course, there is no largest natural number. By short cuts like $10^{10^{10^{...}}}$ we are able to express numbers with precisely known values surpassing any desired magnitude. But it will never be possible, by any means of future technology and of mathematical techniques, to know that natural number *P* which consists of the first $10^{100}$ decimal digits of π (given π is a normal irrational without any pattern appearing in its *n*-adic expansion). This *P* will never be fully available. But what is a natural number the digits of which will never be known? *P* is an idea but it is not a number and, contrary to any 80-digits string, it will never be a number. It is even impossible to distinguish it from that number *P'*, which is created by replacing the last digit of *P* by, say, 5. It will probably never be possible to decide, which of the following relations holds

>*P* < *P'* or *P* = *P'* or *P* > *P'*.

If none of these relations can ever be verified, we can conclude, adopting a realistic philosophical position, that none of them is true.

A recently devised method to compute hexadecimal digits of π without knowing the previous digits [1] fails in this domain too because of finite precision and lacking memory space. But even if the last digit of *P* could be computed, we would need to know all the other digits to distinguish *P* from all similar numbers 314...*d*... where *d* means the *n*-th digit with $1 \leq n \leq 10^{100}$. And why should we stop at the comparatively small Googol?

It is obvious how to apply the aforementioned ideas to the collection of rational numbers. Rational numbers can be defined as equivalence classes of pairs of natural numbers. By multiplying a rational number by its denominator, we obtain a natural number. Natural numbers measure values based upon the unit, rational numbers measure values depending on their denominator. All rational numbers have some terminating *n*-adic expansions. In other representations they have periods. In case the terminating expansion or the period is not too long, these numbers can be identified and hence they do exist. A rational number which approximates π better than $1/2^{10^{100}}$, i.e., to $10^{100}$ bits, does not exist.

## ARE THE REAL NUMBERS REALLY REAL?

2500 years ago Hippasos of Metapont discovered that there are properties which cannot be expressed by ratios of natural numbers. But does this automatically imply that irrational numbers exist? Shouldn't we be able to acknowledge that some things cannot be expressed by numbers? Who would be surprised to learn that some emotions do not fit into the set of numbers? Who would maintain that "justice" is the number 4 as the Pythagoreans insisted?

The task to find the ratio between circumference and diameter of the ideal circle does exist as an idea, and it has a name. Euler called it *p* or *c*, until π became generally accepted - not least supported by Euler's adaption of Jones' π in his *Introductio in Analaysin Infinitorum*. But is this idea, named π, a number? Does the ratio of circumference and diameter exist in reality or in the brains of mathematicians? Surely not. In reality an ideal circle cannot exist because any real circle requires more mass than can be tolerated by an Euclidean space. On the other hand nobody will be able to imagine it. Not even the sharpest mind will be able to determine π better than up to two digits by pure imagination, few will be able to compute more than 5 digits of π in their mind using any desired algorithm, and even fewer will be able to memorize more than the first 50 digits (although it has been shown possible by Hideaki Tomoyori to know by heart 40000 digits) of π.

Cauchy, Weierstrass, Cantor, and Dedekind attempted to give meaning to irrational numbers, well aware that there was more to do than to find suitable names. According to Cantor, √3 is not a number but only a symbol: "√3 ist also nur ein *Zeichen* für eine Zahl, welche erst noch gefunden werden soll, nicht aber deren Definition. Letztere wird jedoch in meiner Weise etwa durch (1.7, 1.73, 1.732, ...) befriedigend gegeben" [2].

It is argued that √3 does exist, because *it can be approximated to any desired precision* by some sequence $(a_n)$ such that for any positive ε we can find a natural number $n_0$ such that $|a_n - \sqrt{3}| < \varepsilon$ for $n \geq n_0$. It is clear that Cantor and his contemporaries could not perceive the *principle* limits of their approach. But every present-day scientist should know that it is *in principle* impossible to approximate any irrational like √3 by $\varepsilon < 1/2^{10^{100}}$. (An exception are irrationals like $\sqrt{3}/2^{10^{100}}$. But that does not establish their existence.) Therefore, the request to achieve *any desired precision* fails in decimal and binary and any other fixed *n*-adic representation.

What about continued fractions, sequences, series, or modular identities determining irrational numbers *x* with "arbitrary precision"? All these methods devised to compare *x* with 1 must necessarily fail, because the *uninterrupted* sequence of natural numbers up to "arbitrary magnitude" required to calculate the terms and to store the rational approximations is not available. As a result we find that irrational numbers do not exist other than as ideas. Irrational numbers simply are not numbers.

**INFINITY FINISHED.**

Transfinite set theory is based upon the observation that there are different infinities. Within its framework there arise quite a lot of paradoxes, for instance the following: By his famous diagonal argument Cantor [3] proved that the set of real numbers (usually only the real interval between 0 and 1 is considered) is non-denumerable. Its cardinal number is $2^{\aleph_0}$ which is distinctly larger than that of the natural numbers, $\aleph_0$. If we consider the binary tree, however,

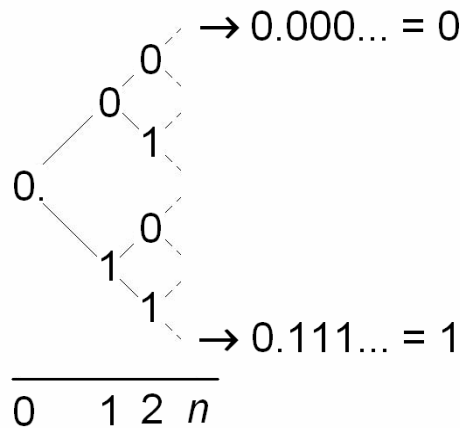

Fig. 1: The binary tree of all real numbers of [0,1].

we find that its nodes, the bits, are countable. Every real number of the interval, including all irrationals and all non-terminating rationals like $1/3 = 0.010101...$, as well as all terminating rationals like $0.1000...$, is realized by one or two infinite paths of the binary tree, each of which begins at "0.", the first node, and does never end. Up to level *n* the tree contains $B(n) = 2^{n+1} - 1$ nodes or bits. $2^n$ separate paths arrive at level *n* while $P(n) = 2^{n+1}$ separate paths leave it. Therefore we have always, i.e., up to any level enumerated by a finite natural number *n*,

$P(n) \leq B(n) + 1$.

But even "in the infinite", should it exist, a path cannot split into two paths without creating a node. Because a node is *defined* as a branching point, no increase in $P(n)$ is possible without the same increase in *B*(*n*). Therefore, the number *P*(*n*) of paths cannot be uncountable unless the number *B*(*n*) of nodes is uncountable too.

We can even establish a bijection between nodes and paths, for instance as follows: That path leaving a node with negative slope is defined to be the continuation of the incoming path whereas the node is mapped on that path leaving it with positive slope. In this way every path is counted as soon as it separates from other paths. We cannot, of course, determine that node which is mapped on √3. But if √3 is an individual number which can appear in a list of real numbers (like those being subject to Cantor's diagonal argument) then it must have its individual path in the tree and its individual node with no doubt.

Being sure that the number *B*(*n*) of nodes is countable the number *P*(*n*) of path must be countable too. By herewith contradicting the result of Cantor's diagonal proof we see that it is not possible to argue consistently in the domain of infinity. Two reasonable proofs can lead to opposite results.

This and other paradoxes are solved by the fact that infinite sets do not exist at all. There is no way to circumvent this conclusion. What the mathematicians of the 19$^{th}$ century could not yet know, and what those of the 20$^{th}$ century seem to have pushed out of their minds: The universe contains less than $10^{80}$ protons and certainly less than $10^{100}$ particles which can store bits. It is, however, a prerequisite of set theory that an element of a set must differ from any other element of that set by at least one property. For this sake one would need at least one bit per element. Therefore, an upper limit of the number of elements of all sets is $10^{100}$. The supremum is certainly less. We have to revise the idea of actually existing Cantorian sets. Even the smallest one, the set ℕ of *all* natural numbers does not exist, let alone the set of all real numbers. At least two notions are to be distinguished in *mathe-realism* with respect to natural numbers.

1) Only those natural numbers exist which are available, i.e., which can be used by someone for calculating purposes. This proposition, including the "someone", is left somewhat vague on purpose. The existence of a natural number has "relativistic" aspects: The question of its existence can be answered differently by different individuals and at different times. As an example consider a poem which exists for the poet who just has written it but not yet for anybody else. Another example is the set of prime numbers to be discovered in the year 2010. It does not yet exist. It is unknown how many elements it will have, even whether it will have elements at all.

The set ℕ* of all natural numbers which exist relative to an observer can increase or decrease like the set of all known primes. Therefore it is difficult to determine its cardinal number |ℕ*|. But obviously |ℕ*| < $10^{100}$ is not infinite.

2) All natural numbers which have ever existed and which ever will exist do not form a set in the sense of set theory because not all of them are simultaneously available and distinguishable. Some of them do not exist yet, like some of the 80-letters poems. The number of all those numbers of this collection is *potentially* infinite because it is not bounded by any threshold, but in no case any cardinal number can become the actual infinity $\aleph_0$ let alone exceed it.

We can conclude: The set $\mathbb{N} = \{1, 2, 3, ... \}$ of *all* natural numbers in the sense of Cantorian set theory does not exist. It is simply not available. But it is difficult to recognize that, because usually only some comparatively small numbers are chosen as examples and in calculations, followed by the dubious symbol "...". The natural numbers, in general thoughtlessly imagined as an unbroken sequence, do not come along like the shiny wagons of a long train. Their sequence has gaps.

As a result, we find in *mathe-realism* that infinite sets consisting of distinct elements cannot exist, neither in the brain of any intelligent being nor elsewhere in the universe. $10^{100}$ is an upper limit for any cardinal number. Therefore, all the paradoxes of set theory, not easy to be circumvented otherwise, vanish.

We give the last word to Abraham Robinson who was a pupil of A.A. Fraenkel, one of the most distinguished set theorists. Nevertheless Robinson recognized already 40 years ago: "(i) Infinite totalities do not exist in any sense of the word (i.e., either really or ideally). More precisely, any mention, or purported mention, of infinite totalities is, literally, *meaningless*. (ii) Nevertheless, we should continue the business of mathematics 'as usual', i.e., we should act *as if* infinite totalities really existed."

**APPENDIX**

The bijection of paths and nodes of the binary tree proposed in this paper

> That path leaving a node with negative slope is defined to be the continuation of the incoming path whereas the node is mapped on that path leaving it with positive slope.

has been accused to cover only those paths which represent rational numbers (ending with an infinite sequence of ones). But obviously no path can split into two paths other than at a node. To assert the existence of more path than nodes is to assert the existence of more split positions than nodes - and that is equivalent to claim the existence of more nodes than nodes. Therefore, if the objection is justified, then it is only because there are no other than rational numbers and, hence, there are no other than the corresponding paths. Nevertheless this objection can be met either by using a random mapping or by a fractional allocation of nodes.

In the first case we know that there is one more node with one more split (because a split position is a node). Therefore it does not matter which path is declared to be that path carrying a node already and which path is in need of a node.

In the latter case, the first node 0 on level 0 is divided into two equal shares. One share is inherited by all the path beginning with 0.0 while the other share is inherited by all the paths beginning with 0.1. Similarly, half of the first node 0 on level 1 is inherited by all the paths beginning with 0.00 while the other half is inherited by all the paths beginning with 0.01. Half of the first node 1 on level 1 is inherited by all the paths beginning with 0.10 while the other half is inherited by all the paths beginning with 0.11. Continuation of this process leads to the following sequence of shares of nodes which are mapped on the finite segments of paths:

> 1/2 node up to level 1
>
> 1/2 + 1/4 nodes up to level 2
>
> 1/2 + 1/4 + 1/8 nodes up to level 3

briefly

> $1 - 1/2^n$ nodes up to level $n$

Continuing this calculation we obtain for the complete infinite path covering every level 1/2 + 1/4 + 1/8 + ... = 1 node.

As the number of nodes is countable, the number of paths is countable too. But the number of paths is not less than the number of real numbers in the interval [0, 1]. This is a contradiction to Cantor's theorem.